\documentclass[a4paper]{elsarticle}

\usepackage{bm}
\usepackage{float}
\usepackage{color}
\usepackage{multirow}
\usepackage{graphicx}
\usepackage{easybmat}
\usepackage{amsmath,amssymb,amsfonts}
\usepackage[margin=3cm]{geometry}
\usepackage{subfig}
\usepackage[table]{xcolor}
\usepackage{soul}

\usepackage{gamosinos}

\definecolor{lightgray}{gray}{0.9}
\usepackage[colorlinks,
            bookmarksopen,
            bookmarksnumbered,
            citecolor=blue,
            linkcolor=blue,
            urlcolor=black]{hyperref}

\nopreprintlinetrue

\usepackage{lineno}

\newcommand{\K}{\bm{K}}
\renewcommand{\u}{\bm{u}}
\newcommand{\f}{\bm{f}}

\begin{document}

\title{On the data-driven description of lattice materials mechanics}

\author[1]{Ismael Ben-Yelun\corref{cor1}}\ead{i.binsenser@upm.es}
\cortext[cor1]{Corresponding author}

\author[1,2]{Luis Irastorza-Valera}

\author[1,3,4]{Luis Saucedo-Mora}

\author[1,5]{Francisco Javier Montáns}

\author[2]{Francisco Chinesta}

\affiliation[1]{organization={E.T.S. de Ingeniería Aeronáutica y del Espacio, Universidad Politécnica de Madrid},
addressline={Pza. Cardenal Cisneros 3},%
city={Madrid},%
postcode={28040},%
country={Spain}%
}

\affiliation[2]{organization={PIMM Laboratory, Arts et Métiers Institute of Technology},
addressline={151 Bd de l'Hôpital},%
city={Paris},%
postcode={75013},%
country={France}%
}

\affiliation[3]{organization={Department of Materials, University of Oxford},
addressline={Parks Road},
city={Oxford},
postcode={OX1 3PJ},
country={UK}
}

\affiliation[4]{organization={Department of Nuclear Science and Engineering, Massachusetts Institute of Technology},
state={Massachusetts},
postcode={MA02139},
country={USA}
}

\affiliation[5]{organization={Department of Mechanical and Aerospace Engineering, Herbert Wertheim College of Engineering, University of Florida},
state={Florida},
postcode={FL32611},
country={USA}
}

\begin{abstract}
In the emerging field of mechanical metamaterials, using periodic lattice structures as a primary ingredient is relatively frequent. However, the choice of aperiodic lattices in these structures presents unique advantages regarding failure, e.g., buckling or fracture, because avoiding repeated patterns prevents global failures, with local failures occurring in turn that can beneficially delay structural collapse. Therefore, it is expedient to develop models for computing efficiently the effective mechanical properties in lattices from different general features while addressing the challenge of presenting topologies (or graphs) of different sizes. In this paper, we develop a deep learning model to predict energetically-equivalent mechanical properties of linear elastic lattices effectively. Considering the lattice as a graph and defining material and geometrical features on such, we show that Graph Neural Networks provide more accurate predictions than a dense, fully connected strategy, thanks to the geometrically induced bias through graph representation, closer to the underlying equilibrium laws from mechanics solved in the direct problem. Leveraging the efficient forward-evaluation of a vast number of lattices using this surrogate enables the inverse problem, i.e., to obtain a structure having prescribed specific behavior, which is ultimately suitable for multiscale structural optimization problems.
\end{abstract}

\maketitle


\section{Introduction}

\subsection{Motivation}
Lattice structures in the form of micro-architected lattices constitute a building block for the emerging field of mechanical metamaterials, which have gained attention recently due to their tunable properties compared to their conventional bulk materials \cite{gibson1999cellular,fleck2010micro,zheng2014ultralight}. This versatility \cite{yang2019novel} in their properties makes them a perfect candidate for structural optimization problems, which poses an additional challenge given that metamaterials span several scales, i.e., yielding the subsequent multiscale problem \cite{ashby1991overview,ashby1993materials}. The application of lattice materials in engineering fields demands real-time calculations, which require the fast evaluation of models, e.g., the mentioned structural optimization \cite{diaz2010topology} or structural health monitoring (SHM) \cite{ozbey2014wireless}. 

In the framework of multiscale optimization problems, different techniques are used for computing the homogenized or equivalent behavior---see Somnic et al. for a more in-depth review on lattice materials \cite{somnic2022homogenization}. To compute the effective behavior, e.g., stiffness, of lattice structures spanning different length scales, numerical techniques addressing the full-resolution of the structure, such as FE$^2$ \cite{schroder2014numerical} emerge to carry out this multiscale problem, but at a high computational cost. Homogenization methods \cite{tessarin2022multiscale} such as FFT homogenization \cite{moulinec1994fast,schneider2021review} overcome this expensive problem. However, this method presents issues when there are two phases with different stiffness, very pronounced in the case of lattice materials since one of the phases is void, i.e., null stiffness---the so-called infinite stiffness contrast in FFT homogenization methods \cite{danesh2023challenges}. 

Data-driven techniques bypass this expensive step, thus allowing for fast evaluations by building less complex surrogate or reduced order models \cite{montans2019data} that learn the (potentially) complex underlying physical phenomena, e.g., using machine learning (ML) or deep learning (DL) algorithms. Once the surrogate models have been trained in an \textit{offline}\ADDED{stage,}they present benefits in different applications in the \textit{online}\ADDED{stage,}e.g., metamaterial characterization and optimization, boosting the computation speed for additional simulations within a multiscale optimization framework, or fast prediction in real-time scenarios, among others. These approaches are widely applied in other engineering fields, like solving complex numerical problems such as parameterized partial differential equations (PDEs) \cite{brunton2016discovering,hernandez2021deep,vinuesa2021high}. 

A popular technique in recent years is the so-called physics-informed neural network (PINN) \cite{raissi2019physics,monaco2023training}, directly imposing the physical equations in the problem, e.g., by introducing the corresponding residuals in the loss function. PINNs ultimately allow for the solution of forward and inverse problems \cite{gao2022physics,di2023physics}---noting the ill-posedness of the latter  \cite{lavrentiev1986ill, kaipio2006statistical}. Initially, the challenges tackled by DL use dense neural networks (DNNs), achieving a good level of precision \cite{kollmann2020deep}. However, they typically require many trainable parameters and a fixed, i.e., dense, input vector. Thus, other architectures such as graph neural networks (GNNs) \cite{sanchez2021gentle}, overcome this step by dealing with variable-size input data (graphs), succeeding to learn from simple essential mechanics such as linear elasticity up to more complex non-linear problems \cite{pichi2024graph}, using a smaller number of parameters.

\ADDED{GNNs have proven useful to tackle graph-like structures, both from a conceptual -- knowledge graphs \cite{Ye2022} describing relationships between nodes e.g., social networks \cite{Fan2019} -- and practical point of view---intertwining nodal and edge behavior as in molecular behavior \cite{Gilmer2017}, neural interaction \cite{Wein2021} or structural mechanics \cite{Parisi2024}. These features can help process what regular DNNs cannot, such as non-Euclidean data presented in a graph form \cite{Wu2021} or tree dependencies not captured in an identically layered network \cite{Zhou2020}. This hardware approach can convey the physical layout of the predicted model more accurately than DNNs -- or any other sort of ML architecture, for that matter -- which is convenient when developing accurate mock representations of the target system i.e., Digital Twins. Diffusion phenomena constitute a clear example of an easily implementable model tailored for graph layouts \cite{Chamberlain2021} with varied applications e.g., heat transfer.}

\subsection{Related research}
Most metamaterial lattices in the literature are periodic and symmetric---built by repetition of a parametric unit cell, which significantly eases their design and manufacturing. Nevertheless, that fact also makes them vulnerable to the propagation of mechanical issues such as buckling, fatigue, or creep. Having identical struts orderly distributed throughout space means that, should any of them reach critical conditions, the rest will contribute to transform a local failure into a global one. 

A solution to that would be an aperiodic lattice so that no two struts would be identical, and thus, never share the exact\ADDED{same}critical requirements. Aperiodic lattices also present advantageous applications in fracture, e.g., toughness improvement \cite{choukir2023role} and vibration suppression \ADDED{via} enhanced damping \cite{liu2022enhanced}. While aperiodic \ADDED{configurations} are being designed and tested \cite{Yang2017, DAlessandro2020}, a unified framework or methodology relating their architecture to their final properties remains elusive, let alone inverse analysis -- obtaining a tailored architecture from the desired performance -- although some promising ML-based techniques have been developed \cite{Bastek2021}.

In periodic metamaterials, the geometry and topology of the unit cell are tightly defined beforehand. However, more rigorous classifications have been proposed \cite{Zok2016} and optimized \cite{Messner2016}. Their properties can usually be obtained -- at least, partially -- from their geometry via experimental testing \cite{Meza2017}, the Finite Element Method (FEM) \cite{Zhang2018} and NNs \cite{Wei2021}. This is different for (semi)randomly-generated architectures, so the first step would be to find a robust way to classify them i.e., telling them apart from each other so their relative advantages and inconveniences can be evaluated. 

Some statistical descriptors are already used for heterogeneous media like alloys or composites \cite{Torquato2002,Cui2021}. Topological Data Analysis (TDA) \cite{Munch2017,Chazal2021}, despite being more widespread in scenarios where one wishes to make sense out of and/or correlate data cloud points to extract the underlying behavior law -- such as contagion scenarios \cite{Taylor2015} -- has been successfully used to classify composites \cite{Runacher2023}. Nonetheless, some studies suggest a topological perspective may not suffice to capture nor predict the metamaterial's mechanical behavior \cite{Meza2017} and call to incorporate homogenization \cite{Vigliotti2014}, affinity deformation considerations between different ``phases'' (unit cells) \cite{Mirzaali2020} or even ML to capture non-linear responses \cite{Xue2023}. Homogenization is ever-important in lattices as a method for accurate, simplified behavior portrayal \cite{Hassani1998-1,Hassani1998-2,Hassani1998-3}, either in a purely analytical form \cite{Coutris2020,Braides2023} or in reverse \cite{Xu2021}, via FEM \cite{Vladulescu2020} or graph-assisted with varied applications \cite{Johnson2018,DeAmbroggio2020,Yang2023,Dold2023}.

ML-based tools such as Convolutional Neural Networks (CNNs), Genetic Algorithms, or DL have been applied in metamaterial design \cite{Ji2022} or fracture mechanism forecast \cite{Lew2021}. In particular, GNNs \cite{Scarselli2009} and Message-Passing \cite{Gilmer2017} have proven very useful for capturing -- and optimizing \cite{Seo2023} -- the topology while learning the underlying physics and the respective conservation laws. They have been used in surrogate modeling \cite{Jiang2023} and prediction of mechanical behavior such as buckling \cite{Prachaseree2022}, fracture \cite{Karapiperis2023,Perera2023}, fatigue \cite{Thomas2023} and non-linear stress/strain relationships \cite{Maurizi2022}, which are crucial needs in SHM.

GNNs can be further enriched by introducing the PDEs describing the physical phenomenon to be replicated in a similar fashion to the aforementioned PINNs \cite{Raissi2019,Karniadakis2021}---which can also be expanded by graph theory \cite{Hall2021}. Some variants of this mixture have tackled unstructured meshes \cite{gao2022physics}, natural convection in fluid dynamics \cite{Peng2023}, microneedle design \cite{Chumpu2023}, soft-tissue mechanics \cite{Dalton2023} or production forecasting \cite{Liu2023}.

All these tools can be viewed as surrogate models since they bypass the need for accurate analytical expressions to describe or predict the material's behavior, mechanical in this case. Delving into non-deterministic problems, they pose an extra issue since the outcome is not known in advance, so models -- such as neural networks -- cannot be ``trained'' as easily. Although bibliography on the topic remains scarce, some non-deterministic/stochastic surrogate examples can be found \cite{Das2021,Jordan2022,Luethen2023}, even for the prediction of material properties \cite{Feng2023,HaririArdebili2023}; namely fatigue \cite{Wang2023}, fracture \cite{Cheng2023} and corrosion \cite{huzni2022physics}.

\subsection{Materials and methods}

In this paper, we build a surrogate mechanical model to predict the equivalent, effective behavior of lattice materials, thus alleviating the computationally expensive associated cost compared to, e.g., an FE$^2$ approach \cite{schroder2014numerical}. To this end, we present and describe the \textit{offline}\ADDED{stage}of this process, showcasing some examples of applications with a synthetically generated dataset. We create our dataset with randomized, unit cell structures due to the highlighted properties that they\ADDED{exhibit.} 

We set a number of nodes\ADDED{randomly distributed}within a cubic domain and connect them through Delaunay triangulation \cite{Bowyer1981, Watson1981}, using either pin-jointed trusses or Euler-Bernoulli beams to model the members/bars of the structures. This way, the lattice structures generated present different numbers of nodes and struts, which poses a challenge for architectures presenting a fixed number of inputs, e.g., regular DNNs. All the structures are labeled with their mechanical equivalent properties, namely Young's moduli $E_z^{\rm eq}$ in $z-$axis, computed using an energy-preserving method as the relationship between the reaction forces and a prescribed displacement.

The effective behavior is computed by fitting an ML surrogate mechanical model, similar to other research conducted on this topic\ADDED{\cite{Kim2021,Shen2022,ben2023gam,Lew2023,Meier2024},}adding the challenge of developing a model able to deal with lattices/meshes of different sizes and resolutions.\ADDED{This paper represents a proof of concept approach in which we are considering linear elastic problems, and then move towards non-linear cases as future steps. An extension to non-linear constitutive behavior could be performed by including more features in the model. For instance, in path-independent non-linear behaviors e.g., hyperelasticity, we hypothesize that the strain energy curve description should be quantified as extra features to include. Likewise, an approach for path-dependent settings e.g., plasticity, would also require the loading history, defining the extra features in terms of internal variables.}Regarding the surrogate model, we propose two non-intrusive, NN architectures to predict the so-mentioned equivalent mechanical property: (1) a DNN approach fed with `engineered' features, and (2) a GNN model using the geometry and material of the structure as inputs.

We first explore the use of a DNN aimed at predicting the mechanical properties with the minimum physical information, e.g., not providing the assembled stiffness matrix but rather the geometry or material properties of the lattice struts. In order to deal with samples presenting different (input) dimensions -- number of nodes, number of members -- reduced techniques shall be applied, e.g., projecting, averaging, interpolating, or pooling.  One option is the proper orthogonal decomposition with artificial NN, i.e., POD-ANN \cite{san2019artificial}, which consists\ADDED{in}projecting the input space in a lower dimension manifold $\mathcal{M}\subset \mathbb{R}^d$, where $d$ is fixed so the latent space of all structures can be fed\ADDED{to}the same dense architecture. This can be achieved using linear, e.g., PCA \cite{bishop2006pattern}, or non-linear, e.g., autoencoder  techniques---similar to \cite{pichi2024graph}. 

This approach becomes a challenging task, leading to the loss of crucial information in the encoding, such as the connectivity of the members within the lattice. Therefore, the DNN approach addressed in this paper is different. A structure is sliced into a certain number of slices in its three directions.\ADDED{For}each slice, the sum of effective areas of intersected bars is computed. The total effective area per slice is fed to a DNN model with $3n_s$ inputs, where $n_s$ is the number of slices selected. Thus, we elaborate a physical analogy with a structure with multiple one-dimensional springs stacked in series. Each slice would represent springs in parallel---hence the sum. An illustrative example of this analogy is described in Section~\ref{subsec:toy_example}. 

This way, we provide the maximum physical knowledge without directly feeding the (assembled) stiffness matrix of the structure,\ADDED{aiming at bringing a physics-informed prior to address the problem.}However, the previous perspective lacks topological information on the structures, hence the exploration of architectures accounting for the connectivity of the nodes as the training model, like GNN. With a graph-level regression GNN, we show an even higher level of accuracy using only as input the coordinates of the nodes, the connectivity, and properties of the members. We note that whereas DNNs explain better what it means to learn the material, GNNs succeed in learning the structure. Taking all the previous into account, we derive a model for computing effective properties, improving the commonly used models that are based on empirical exponential laws relating the stiffness to strut radius, i.e., the relative density of the cell \cite{Meza2017}. In fact, due to the non-bijectivity of the model, the density can be considered an extra degree of freedom.

Once surrogate models of sufficient quality are obtained, the inverse analysis might be performed by efficiently (forward-)evaluating a large number of samples,\ADDED{bypassing the computational expensive step of solving the equilibrium equation.}Furthermore, the non-bijectivity of the problem allows for multiple lattices behaving similarly, thus leaving room to optimize a second variable, such as the weight of the structure, in the subsequent multi-objective optimization problem that might be addressed making use of Pareto optimality \cite{boyd2004convex}. This addresses the issue of multiscale optimization proposed in \cite{ben2023topology}, with the advantage of replacing the macro-to-micro step dependent on empirical laws in the latter with the forward-evaluation of the derived surrogate models proposed in this paper.


\ADDED{As a summary, this paper provides a novel bottom-up methodology i.e., particular to general, from the design of the metamaterial's architecture to the \textit{offline} data-driven description of its properties -- namely its equivalent Young's modulus $E_{z}^{\rm eq}$ -- robust enough to tackle both periodic and aperiodic layouts by means of several ML techniques, comparing their performance and reinforcing the idea of GNNs as the most efficient.}

\ADDED{Whereas other approaches have opted to directly incorporate mechanical and geometrical variables into their ML model e.g., detailed boundary conditions or angle-bending forces \cite{Reid2018}, the presented method yields accurate mechanical behavior predictions even when just provided with very basic geometrical information regardless of the topological complexity or the chosen mechanical framework i.e., pin-jointed nodes, or frames. In doing so, our approach handles aperiodic configurations ignoring the periodic elementary cell \cite{Meier2024}, or any hierarchical orthotropic underlinings, like in composites fibers \cite{Zhu2024}. 
This framework draws inspiration from irregular designs found in nature, e.g. biological designs and heterogeneous microstructures in porous media, which have been addressed by ML techniques but still within periodic unit cell limitations \cite{Shen2022,Nguyen2022}.}

\ADDED{Considering the ill-posedness of the  inverse analysis, we propose a straightforward way to find existing lightweight architectures with a desired property i.e., equivalent Young's modulus, within a close interval from the target, out of a discrete lattices portfolio. Although simple, it may prove a more convenient and versatile method than existing ones---mostly either periodic \cite{Qin2018,Zheng2024} or limited to 2D designs \cite{Zheng2021}. }

\section{Lattice structures database}
\subsection{Structure definition}
Each lattice structure is defined in a unit cube $\Omega \subset \mathbb R^3$ and\ADDED{is}made of $n_j$ nodes (or joints) and $n_e$ members (or struts). Eight nodes of the lattice correspond to a vertex of the cube, thus\ADDED{the only ones with a}fixed position \ADDED{-- unlike \cite{Zheng2023} --} whereas the remaining $n_j - 8$ user-prescribed joints are randomly distributed within the domain $\Omega$. In this paper, samples are generated varying this number from 1 to 50. The 8 fixed nodes act as clamps in which forces or displacements will be applied, being the only nodes whose coordinates are fixed for all experiments so their results can be compared. The rest of the nodes' coordinates are sampled from three independent and identically distributed (iid)\ADDED{uniform}distributions $\mathcal U \left(\ADDED{\delta}, 1-\ADDED{\delta}\right)$, where $\ADDED{\delta}$ is a small value to prevent the randomly located joints from being on the boundary $\partial\Omega$ of the domain.

The elements are defined by a 3D generalization of Delaunay's triangulation (through Bowyer-Watson's algorithm \cite{Bowyer1981,Watson1981}), yielding a total $n_e$ elements, hence this number is not prescribed. Delaunay's method, common in simulation software like ANSYS\copyright, ensures both the absence of unforeseen -- and undesired -- intersections (nodes) and the strut lengths being constrained to a given interval despite their variability. In mechanical terms, it also reduces buckling (struts' slenderness is contained within the aforementioned interval) and prevents the generation of sliver triangles presenting acute angles which would cause stress concentration and instabilities -- mechanisms -- if the truss is pin-jointed.

\begin{figure}[ht]
    \centering
    \subfloat[Imposed displacement and symmetry planes.]{%
    \includegraphics[width=.45\linewidth]{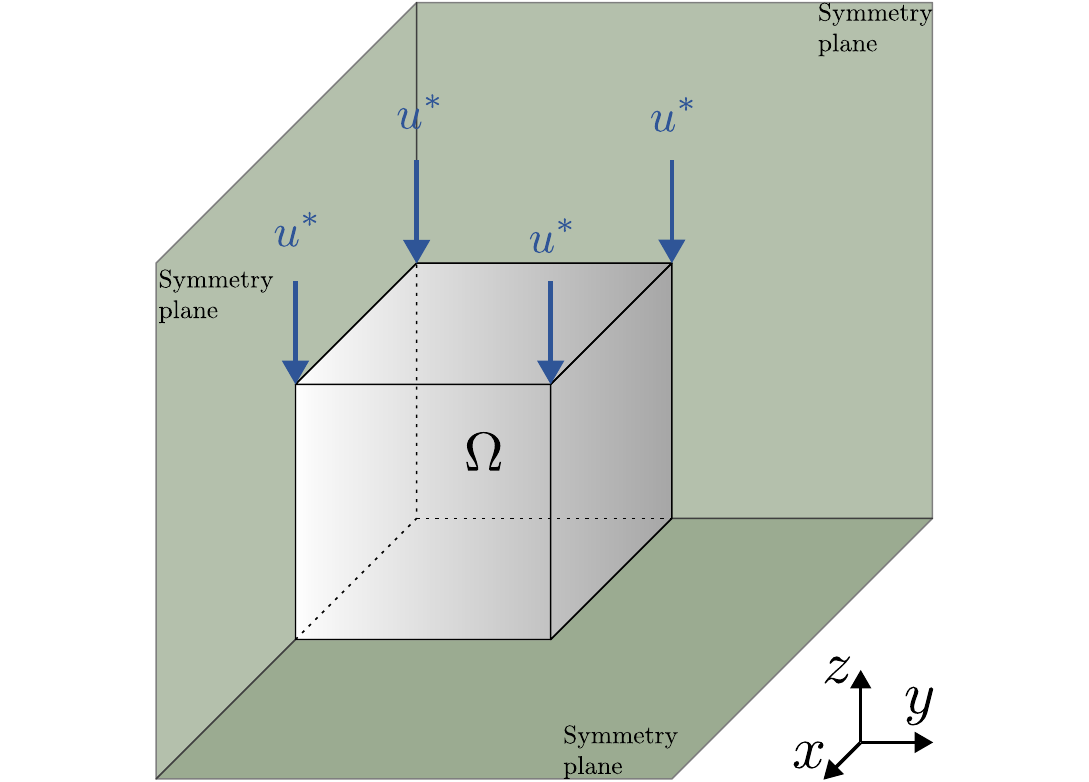}}
    \subfloat[Resulting boundary conditions.]{%
    \includegraphics[width=.45\linewidth]{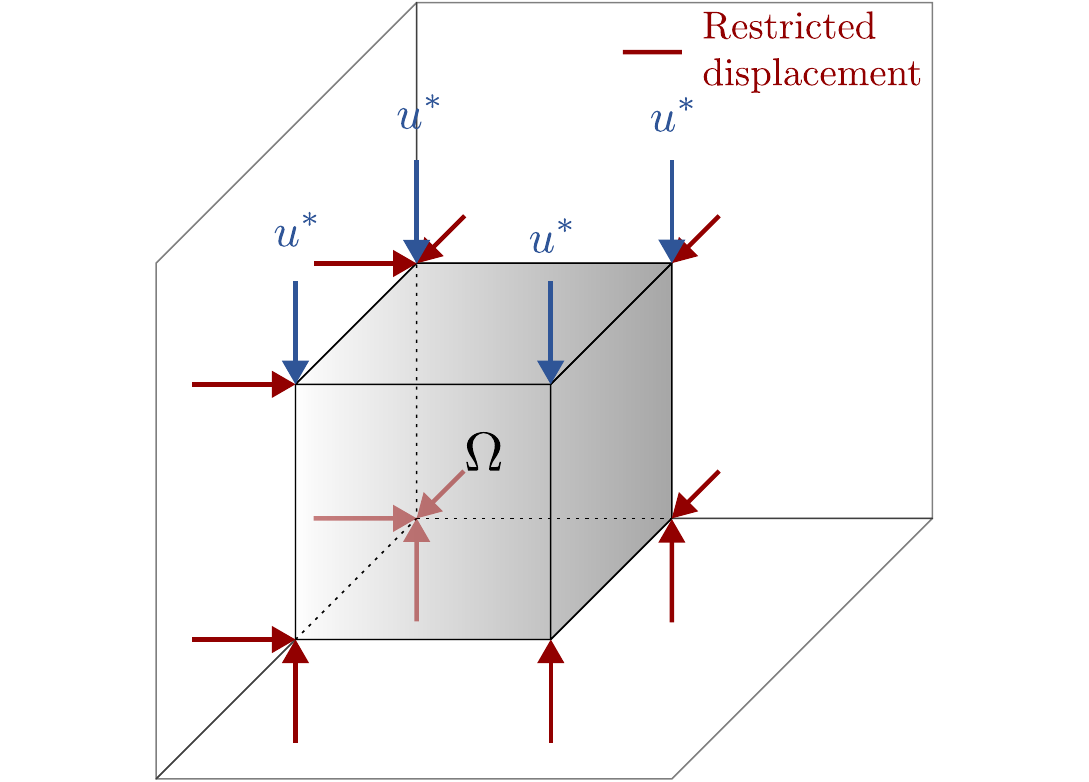}}
    \caption{Boundary conditions of the lattice structure, where $u^*$ is the imposed displacement, in blue. The symmetry planes represented in green in (a) imply the restricted displacements in (b), depicted in red.}
    \label{fig:boundary-conditions}
\end{figure}

\ADDED{\subsection{Boundary conditions and mechanical model of the lattice}}
A lattice is subjected to a uniaxial test with an imposed displacement $u^*$ along the vertical direction ($z-$axis), as it is shown in \figurename~\ref{fig:boundary-conditions}. Due to the three symmetry planes, some degrees of freedom are restricted (or fixed).

In the following, the free degrees of freedom are denoted by $f$, and the restricted ones by $r$. The latter comprise both the imposed and fixed displacements---blue and red arrows in \figurename~\ref{fig:boundary-conditions}, respectively.


The elements of the lattice are modeled either as a pin-jointed truss structure, or an Euler-Bernoulli beam. The (linear) equilibrium equation for a lattice structure with $n_d$ degrees of freedom -- equal to $3n_j$ for trusses, and $6n_j$ for beams --  reads 

\begin{equation}
\K \u = \f,
\label{eq:equilibrium_equation}
\end{equation}
where $\u\in \mathbb R^{n_d}$ and $\f\in \mathbb R^{n_d}$ are respectively the nodal displacements and the nodal forces vector, and $\K\in \mathbb R^{n_d\times n_d}$ is the global stiffness matrix of the whole structure without having imposed the boundary conditions. The local stiffness matrix $\ADDED{\K_e}$ of a truss $\ADDED{e}$ is defined as

\begin{equation}
    \K_e = \dfrac{E_eA_e}{\ell_e}\begin{bmatrix}
        1 & 0 & 0 & -1 & 0 & 0\\
        0 & 0 & 0 & 0 & 0 & 0 \\
        0 & 0 & 0 & 0 & 0 & 0 \\
        -1 & 0 & 0 & 1 & 0 & 0\\
        0 & 0 & 0 & 0 & 0 & 0 \\
        0 & 0 & 0 & 0 & 0 & 0 \\
    \end{bmatrix},
\end{equation}
where $E_e$, $A_e$ and $\ell_e$ are the stiffness, cross-sectional area and length of member $e$; whereas the 12 degrees-of-freedom Euler-Bernoulli beam local stiffness matrix can be found in Equation (5.116) of \cite{przemieniecki1985theory}. Note that the global stiffness matrix used in this work comes from the direct stiffness method (DSM), but it can be extended to other forms of deriving the stiffness matrix e.g., with FE discretization. 

The equilibrium equation in \eqref{eq:equilibrium_equation} can be split\ADDED{into}free $f$ and restricted $r$ degrees of freedom as

\begin{equation}
    \left[
    \begin{BMAT}[2pt,0pt,1.5cm]{c;c}{c;c}
    \K_{ff} & \K_{fr} \\ \K_{rf} & \K_{rr}
    \end{BMAT}
    \right]
     \left[
    \begin{BMAT}[2pt,0pt,1.5cm]{c}{c;c}
    \u_{f} \\ \u_{r}
    \end{BMAT}
    \right]
    =
    \left[
    \begin{BMAT}[2pt,0pt,1.5cm]{c}{c;c}
    \bm{0} \\ \f_r
    \end{BMAT}
    \right],
\end{equation}
where the nodal forces vector in the free degrees of freedoms has been set to $\f_f = \bm{0}$ since the loading is purely displacements imposition. Performing a static condensation, the vector of nodal displacements of the free degrees of freedom $\u_f$ is

\begin{equation}
\u_f = -\K_{ff}^{-1}\K_{fr}\u_r,
\end{equation}
and the nodal forces on the restricted nodes\ADDED{i.e., the reactions,}are

\begin{equation}
\f_r = \left(\K_{rr} - \K_{rf}\K_{ff}^{-1}\K_{fr}\right)\u_r.
\end{equation}

With that, we compute the work of external forces $\mathcal W_{\rm ext}$\ADDED{being introduced into}the lattice

\begin{equation}
    \mathcal{W}_{\rm ext} = \f \cdot \u = \f_r \cdot \u_r.
    \label{eq:work_external_forces}
\end{equation}

\subsection{Effective behavior}
\label{subsec:effective_behavior}
We define an effective material of volume $|\Omega|$ energetically equivalent to the lattice structure. Then, we perform a uniaxial test and equate the energies to compute the equivalent stiffness of the lattice, specifically the Young's modulus along $z-$axis direction---see \figurename~\ref{fig:boundary-conditions}. The potential of the internal forces of the equivalent solid $\Pi_{\rm int}^{\rm eq}$ is

\begin{equation}
    \Pi_{\rm int}^{\rm eq} = \int_{\Omega}\bm{\sigma}^{\rm eq}\left(\bm{\varepsilon}^{\rm eq}(\bm{x})\right):\bm{\varepsilon}^{\rm eq}(\bm{x})\mathrm{d}V,
    \label{eq:potential_internal_equivalent}
\end{equation}
where $\bm{\sigma}^{\rm eq}$ and $\bm{\varepsilon}^{\rm eq}$ are the Cauchy stress and infinitesimal strain second-order tensors of the equivalent body, related through the (linear) constitutive law $\bm{\sigma}^{\rm eq} = \mathbb C^{\rm eq}: \bm{\varepsilon}^{\rm eq}$, being $\mathbb C$ the fourth-order stiffness tensor. Note that any stress-strain work conjugate might be used since the material is working under its linear regime. Particularizing the  constitutive law for a uniaxial test yields $\sigma_{zz}^{\ADDED{\rm eq}} = E_z^{\rm eq}\varepsilon_{zz}^{\rm eq}$, hence the integral in previous Equation~\eqref{eq:potential_internal_equivalent} is straightforward to compute, namely

\begin{equation}
    \Pi_{\rm int}^{\rm eq} = \sigma_{zz}^{\rm eq}\varepsilon_{zz}^{\rm eq}|\Omega|.
\end{equation}

The imposed displacement $u^*$ is the same for both solids, therefore the (continuum-)equivalent strain $\varepsilon_{zz}^{\rm eq}$ is

\begin{equation}
\varepsilon_{zz}^{\rm eq} = \dfrac{u^*}{L_z},
\end{equation}
where $L_z$ is the longitudinal dimension\ADDED{along}the studied direction $z$, i.e., the vertical size of the structure. By energy balance, the potential of external forces and the potential of internal forces are equal, in both lattice and equivalent solids i.e., $\mathcal W_{\rm ext} = \Pi_{\rm int}^{\rm eq}$, yielding

\begin{equation}
\f_r \cdot \u_r = E_z^{\rm eq} \left(\dfrac{u^*}{L_z}\right)^2 |\Omega|. 
\end{equation}

Taking into account that $|\Omega| = A_zL_z$ where $A_z$ is the cross-sectional area of the equivalent solid perpendicular to the studied direction, the equivalent stiffness in $z-$direction $E_z^{\rm eq}$ is obtained

\begin{equation}
    E_z^{\rm eq} = \dfrac{\f_r \cdot \u_r}{u^*A_z}\dfrac{L_z}{u^*}.
    \label{eq:equivalent_stiffness}
\end{equation}

Note that this is equivalent to\ADDED{computing}the ratio of\ADDED{the}resultant of the reaction forces $\f_r$ at the top face and the cross-sectional area perpendicular to the loading application $A_z$ i.e., the equivalent stress $\sigma_{zz}$, and dividing it by the imposed strain $\varepsilon_{zz}$. 



\section{Dense Neural Networks model}
We now turn to the DL prediction of the equivalent stiffness $E_z^{\rm eq}$, object of study of this paper. By making use of a toy example of linear springs arranged in series, we train the simplest DNN model with the most meaningful variables. Taking this example into account, we propose a procedure based on slicing and weighing a given pin-jointed truss lattice structure in order to predict the mentioned mechanical property by training another DNN model.

\subsection{Prediction of springs arranged in series: a toy example}\label{subsec:toy_example}
\begin{figure}[htb]
    \centering
    \includegraphics[width=.49\linewidth]{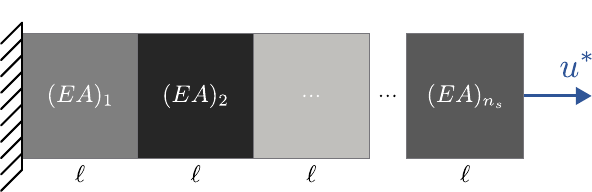}
    \caption{1D linear, elastic problem comprised of $\ADDED{e=1,...,}n_s$ elements of same length $\ell$ and different axial stiffness $(EA)_{e}$ arranged in series.}
    \label{fig:toy_example}
\end{figure}

We consider the following 1D linear, elastic problem comprised of $\ADDED{e=1,...,} n_s$ uniaxial elements with different axial stiffness $(EA)_e$ and equal length $\ell$ arranged in series, as depicted in \figurename~\ref{fig:toy_example}.

By elemental mechanics, the equivalent stiffness $\ADDED{k}_{\rm eq}$ of $n_s$ springs\ADDED{with constant $k_e$}arranged in series is given by

\begin{equation}
    \ADDED{k}_{\rm eq}^{-1} = \sum_{e=1}^{n_s}\ADDED{k}_e^{-1}.
\end{equation}

Considering the relation $\ADDED{k_e = (}EA\ADDED{)_e}/\ell\ADDED{_e}$ between a mechanical spring $\ADDED{e}$ with constant $\ADDED{k_e}$ and a uniaxial element (i.e., bar or truss) with axial stiffness $(EA)\ADDED{_e}$ and length $\ell\ADDED{_e}$ -- which is also known from the theory of pin-jointed truss lattices -- the equivalent\ADDED{axial}stiffness $\ADDED{(EA)_{\rm eq}}$ of the system in \figurename~\ref{fig:toy_example} is

\begin{equation}
    (EA)_{\rm eq} = \dfrac{n_s}{\sum_e (EA)_e^{-1}}.
    \label{eq:spring_series_toy}
\end{equation}

\ADDED{Recall that $(EA)_e$ is the axial stiffness of element $e$.}Therefore, this is the formula that the model would have to learn in case of approaching this problem with an\ADDED{ML approach.}By means of the universal approximation theorem \cite{hornik1989multilayer}, this formula may be accurately represented by a sufficiently wide DNN architecture when given the\ADDED{axial}stiffnesses of the members $(EA)_e$ as inputs. Thus, a DNN model is trained to prove the previous statement. The input features of this model are the inverse\ADDED{axial}stiffness $1/(EA)_e$ of the elements since this allows the DNN to better approximate the target to predict $(EA)_{\rm eq}$ i.e., the function in Equation~\eqref{eq:spring_series_toy}. The architecture used is described in \tablename~\ref{tab:dnn_architecture_toy}. In the following, we consider $n_s = 19$.

\begin{table}[htb]
    \centering
    \begin{tabular}{c|c c}
        Layer & Neurons & Activation \\ \hline
        Input layer & $n_s = 19$ & selu \cite{klambauer2017self}  \\
        \rowcolor{lightgray}Hidden layer \#1 & $5n_s = 95$ & selu \cite{klambauer2017self} \\
        Output layer & $1$ & linear 
    \end{tabular}
    \caption{Toy example. Hyper-parameters of the fully-connected DNN architecture.}
    \label{tab:dnn_architecture_toy}
\end{table}

Similarly, the training (hyper-)parameters are displayed in \tablename~\ref{tab:dnn_optimization_parameters_toy}.\ADDED{We select the nadam optimizer \cite{dozat2016incorporating}, although other available choices such as adam, RMSprop, or adagrad perform similarly in the models proposed. Therefore, we keep nadam as the optimizer for the models trained in this paper.}

\begin{table}[htb]
    \centering
    \begin{tabular}{c|c}
         Parameter & Choice \\ \hline
         Optimizer & nadam \\
         \rowcolor{lightgray}Learning rate & $0.001$\\
         Epochs & $2\,000$ \\
         \rowcolor{lightgray}Loss, $\mathcal L$ & MSE \\
         Validation split & 15\%
    \end{tabular}
    \caption{Toy example. (Hyper-)parameters of the training.}
    \label{tab:dnn_optimization_parameters_toy}
\end{table}

Furthermore, a learning rate reduction on plateau by $0.5$ is applied, as well as an early stopping criterion, both on the loss function. Lastly, the dataset is generated by drawing $N = 20\,000\times n_s$ samples from a\ADDED{uniform}distribution $p(EA) \sim \mathcal U ((EA)_{\min}, (EA)_{\max})$, where $(EA)_{\min}=10$\!\ADDED{N} and $(EA)_{\max}=100$\!\ADDED{N.}This way, the following set is generated 
\[
\ADDED{\left\lbrace\lbrace(EA)_e^{(i)}\rbrace_{e=1}^{n_s},(EA)_{\rm eq}^{(i)}\right\rbrace_{i=1}^{N}.}
\]

A random train-test split of 15\% is then applied, and both the features and targets are normalized with a standard scaler. 

\begin{figure}[htb]
    \centering
    \includegraphics[width=.4\linewidth]{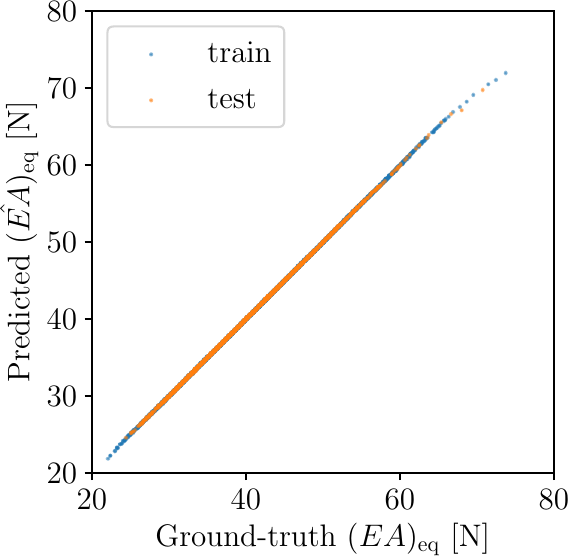}
    \caption{Toy example. Equivalent\ADDED{axial}stiffness $\ADDED{(}E\ADDED{A)}_{\rm eq}$ prediction, representing train and test in blue and orange, respectively. Although there is a slight under-prediction close to the upper limit, the precision is virtually perfect, with coefficients of determination $R^2_{\rm train}=\ADDED{1.0000}$ and $R^2_{\rm test}=\ADDED{0.9999}$.}
    \label{fig:results_toy_example}
\end{figure}

Once the training is performed, both the train and test metrics are computed. Denoting the prediction with a hat $\hat \cdot$,\ADDED{the scaling with a tilde $\tilde \cdot$, and defining the error of a sample $i$ as the squared residual $r^{2}_{i}$, namely}

\begin{equation}
    \ADDED{r^{2}_i = \left({\tilde{E} {}_z^{\rm eq}}^{(i)} - \widehat{{\tilde{E} {}_z^{\rm eq}}^{(i)}}\right)^2},
    \label{eq:residual}
\end{equation}
\ADDED{we display the minimum and maximum errors in addition to the mean error i.e., the loss function $\mathcal L$ or MSE, to assess every model shown.}The DNN yields losses of $\mathcal L_{\rm train}=\ADDED{3.94\cdot 10^{-5}}$ and $\mathcal L_{\rm test}=\ADDED{4.89\cdot 10^{-5}}$, and coefficients of determination of $R^2_{\rm train}=\ADDED{1.0000}$ and $R^2_{\rm test}=\ADDED{0.9999}$.\ADDED{The maximum and minimum errors as defined in \eqref{eq:residual} are $(r_{\min}^2,r_{\max}^2)_{\rm train} = (2.87\cdot 10^{-14}, 6.64\cdot 10^{-2})$ and $(r_{\min}^2,r_{\max}^2)_{\rm test} = (1.34\cdot 10^{-13}, 1.94\cdot 10^{-2})$.}The train and test predictions are shown in a scatter plot in \figurename~\ref{fig:results_toy_example}, where the perfect prediction is given by the line at $45^{\circ}\!$\ADDED{, noting a slight under-prediction close to the upper limit.}

Of course, this is a simple problem solved with an explicit formula, which is virtually perfectly approximated with a DNN model. We use the learnings from this method to address the effective behavior  with a similar surrogate DNN model.

\subsection{Slice features of the lattice structure}
\begin{figure}[htb]
    \centering
    \includegraphics[width=.9\linewidth]{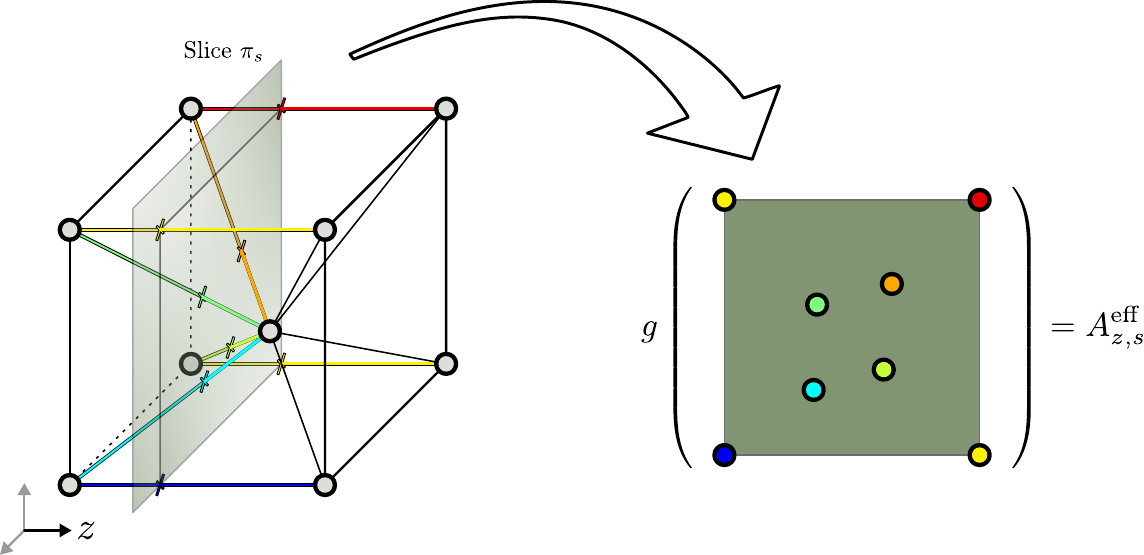}
    \caption{\ADDED{Schematic}drawing of the lattice slicing perpendicular to $z-$axis for an arbitrary slice plane $\pi_s$. The joints are represented in gray, and the intersected bars are represented using a jet colormap. The effective area $A_{z,s}^{\rm eff}$ of the slice $s$ transverse to $z$ is computed by applying a function $g$ yet to be defined.}
    \label{fig:slice-lattice}
\end{figure}

Analogously to the simple idea of the previous toy example, the (random) lattice structure of a given material is now sliced into $n_s$ parts. We formulate the following hypothesis: by defining an \textit{effective} area per slice ($n_s$ in total), the prediction can be performed through an analogous approach to the toy example i.e., a DNN whose inputs are the inverse values of the effective areas. A schematic drawing of a general slice $s$ is depicted in \figurename~\ref{fig:slice-lattice}.

To define the effective area per slice, we introduce the function $g(\cdot)$, which is applied to the areas and directions of the intersected bars by the slice. Given an element $e$ with cross-sectional area $A_e$, and defined with its unit vector by $\bm n_e = [\cos \theta_{x,e},\cos\theta_{y,e}, \cos \theta_{z,e}]^\intercal$, the effective area $A_{z,s}^{\rm eff}$ of the slice $s$ perpendicular to $z$ direction is computed as the sum of the projected areas along $z$ of the intersected bars in the plane $\pi_s$ i.e.,

\begin{equation}
    A_{z,s}^{\rm eff} = \sum_{e\in \pi_s} A_e \cos^2 \theta_{z,e}.
    \label{eq:effective_area}
\end{equation}

\begin{figure}[htb]
    \centering
    \includegraphics[width=.35\linewidth]{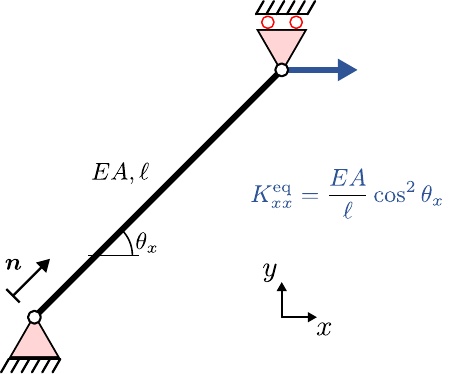}
    \caption{Equivalent stiffness $K_{xx}^{\rm eq}$ along direction $x$ in an example strut to show the quadratic exponent in the cosine function for computing the effective area. The bar is made of a material with Young's modulus $E$, has cross-sectional area $A$ and length $\ell$, and presents an arbitrary inclination with respect to $x$ defined by the angle $\theta_x$.}
    \label{fig:equivalent_stiffness}
\end{figure}

The quadratic exponent in the cosine function in Equation~\eqref{eq:effective_area} comes from the projection of both the forces and displacements onto the local axis of the bar. An illustration of this is depicted in \figurename~\ref{fig:equivalent_stiffness}. By simple rotations, the equivalent stiffness $K_{xx}^{\rm eq}$ in $x$ direction is obtained as the actual stiffness of the bar i.e., $EA/\ell$, multiplied by the cosine square of the angle between the bar and $x-$axis, namely

\[
K_{xx}^{\rm eq} = \frac{EA}{\ell}\cos^2\theta_x.
\]

The same definition applies to any other direction $d$, setting as input the corresponding cosine of the angle $\cos\theta_{d}$. Lastly, the transverse effective areas have to be considered when generalizing the problem from 1D to 2D or 3D. 

\begin{figure}
    \centering
    \includegraphics[width=.7\linewidth]{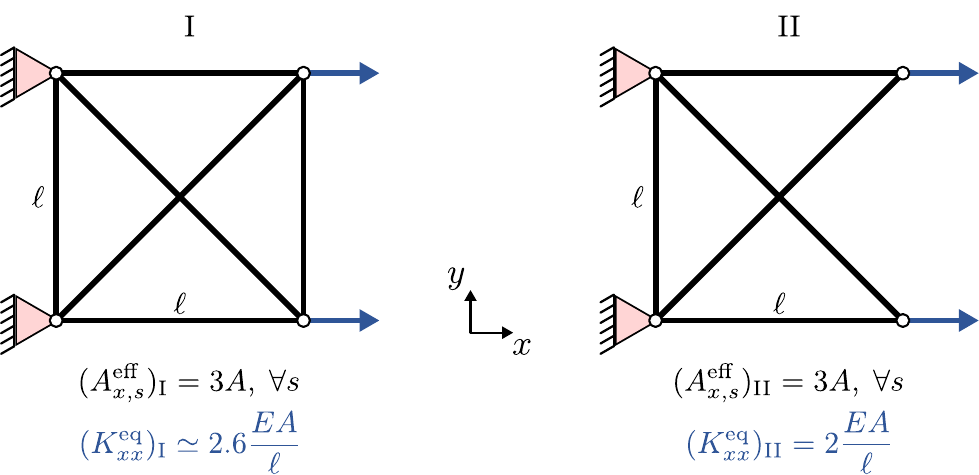}
    \caption{2D lattice structures I and II made of different (square) unit cells of length $\ell$. Whereas the effective area transverse to $x$ is the same, $(A_{x,s}^{\rm eff})_{\rm I} = (A_{x,s}^{\rm eff})_{\rm II} = 3A$ (for every slice $s$), the equivalent stiffness differs, $(K_{xx}^{\rm eq})_{\rm I} \simeq 2.6EA/\ell$ and $(K_{xx}^{\rm eq})_{\rm II} = 2EA/\ell$. All bars are of the same material with Young's modulus $E$ and same cross-sectional area $A$.}
    \label{fig:transverse_effect}
\end{figure}

Let I and II be two 2D lattice structures, with slightly different (square) unit cells of length $\ell$, depicted in \figurename~\ref{fig:transverse_effect}. All bars are of the same material with Young's modulus $E$ and same cross-sectional area $A$. Both lattices I and II have the same effective area perpendicular to $x$ i.e., $A_{x,s}^{\rm eff} = 3A$ for every slice $s$ contained in the structure, according to Equation~\eqref{eq:effective_area}. However, the equivalent stiffnesses $K_{xx}^{\rm eq}$ following the procedure derived in Section~\ref{subsec:effective_behavior}, are different in cases I and II i.e., $(K_{xx}^{\rm eq})_{\rm I} \neq (K_{xx}^{\rm eq})_{\rm II}$, due to transverse effects---lattice II is the same as lattice I with a (transverse) bar removed. Therefore, the effective areas $A_{y,s}^{\rm eff}$ perpendicular to $y$ have to be considered likewise in the surrogate model to make effective predictions for the general, non-1D case. Note that these effective areas for the lattices in \figurename~\ref{fig:transverse_effect} are $(A_{y,s}^{\rm eff})_{\rm I} = 2A$ and $(A_{y,s}^{\rm eff})_{\rm II} = A$.

\subsection{Results}
Using the inverse of the effective transverse areas in the three spatial directions $A_{x,s}^{\rm eff}$, $A_{y,s}^{\rm eff}$, and $A_{z,s}^{\rm eff}$ for $n_s$ slices as inputs, a DNN model for 3D pin-jointed truss lattices is generated analogously to the toy example case shown in Section~\ref{subsec:toy_example}. The dataset is comprised of $10\,000$ randomly-generated samples i.e., lattices. The base material is AISI Type 316L Stainless Steel, annealed bar, with a Young's modulus of $E_{\rm base} = 193\,$GPa, and the lattice members are circular bars with radius $r=5\,$mm, leading to a cross-sectional area of $A = 78.54\,$mm$^2$.

In order to select a suitable number $n_s$ of slices, the following limit analysis is performed on this number:\ADDED{a pseudo-\!\!}equivalent stiffness of the lattices of the dataset is computed via Equation~\eqref{eq:spring_series_toy}, using the effective area transverse to $z$ per slice $s$, $A_{z,s}^{\rm eff}$. Recall that we do not state that this is the actual equivalent stiffness to predict, but serves as a global value to make comparisons for this limit analysis. Then, these samples are sliced using different numbers of slices $n_s$. 

\begin{figure}
    \centering
    \includegraphics[width=.55\linewidth]{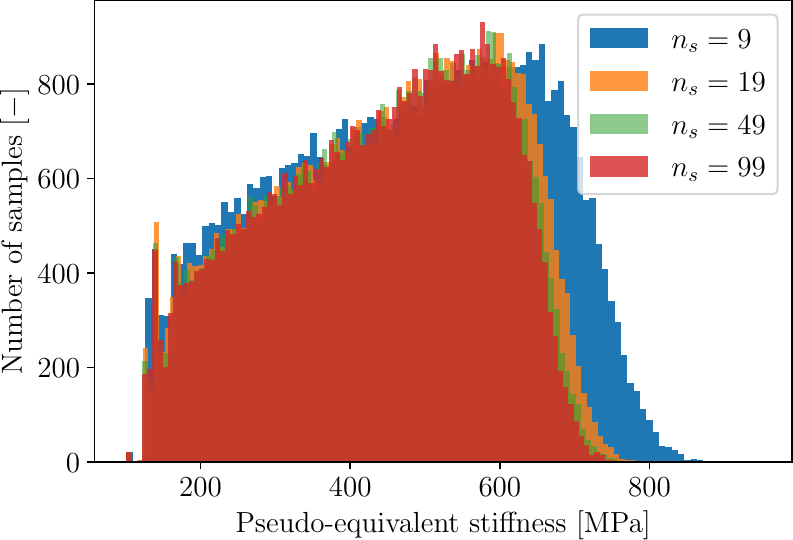}
    \caption{Pin-jointed truss lattice DNN model. Distribution of the\ADDED{pseudo-\!\!}equivalent stiffness according to Equation~\eqref{eq:spring_series_toy} for different number of slices $n_s = \lbrace 9,19,49,99\rbrace$.}
    \label{fig:limit_analysis}
\end{figure}

The histogram of the\ADDED{pseudo-\!\!}equivalent stiffness for the different $n_s$ used is depicted in \figurename~\ref{fig:limit_analysis}. Note that from $n_s=19$ on, this global magnitude per lattice is similarly distributed.\ADDED{Therefore, we train one model considering $n_s=19$ slices and another one with $n_s=49$ to show the performance balance between both.}

Lastly, the equivalent stiffness to predict, $E_{z}^{\rm eq}$ -- Equation~\eqref{eq:equivalent_stiffness} -- is divided by the volume of the lattice to facilitate the predictions of the model. Namely, for each sample

\begin{equation}
    \overline{E_{z}^{\rm eq}} = \dfrac{E_{z}^{\rm eq}}{V},
\end{equation}
where $V$ is the volume of the lattice i.e., 
\begin{equation}
    V = \sum_e A_e \ell_e,
    \label{eq:lattice_volume}
\end{equation}
leading to the following dataset  

$$\left\lbrace \lbrace (A_{x,s}^{\rm eff})^{(i)}, (A_{y,s}^{\rm eff})^{(i)}, (A_{z,s}^{\rm eff})^{(i)}\rbrace_{s=1}^{n_s},\overline{E_{z}^{\rm eq}}^{(i)}\right\rbrace_{i=1}^{N}.$$

\begin{figure}
    \centering
    \includegraphics[width=.4\linewidth]{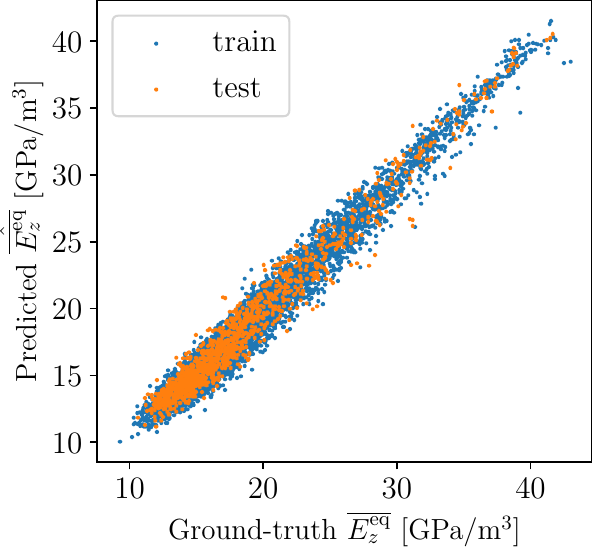}
    \caption{Pin-jointed truss lattice DNN model. Equivalent stiffness per unit volume $\overline{E_{z}^{\rm eq}}$ prediction, where the train and test are represented by blue and orange points, respectively.}
    \label{fig:results_dnn_example}
\end{figure}

The same architecture and optimization parameters displayed in Tables \ref{tab:dnn_architecture_toy} and \ref{tab:dnn_optimization_parameters_toy} are used. A validation set of 15\% of the samples is separated, applying a train/test 85-15\% split in the remaining dataset. Both the features and targets are normalized with a standard scaler. Once the training is performed, both the train and test metrics are computed. The DNN yields losses of $\mathcal L_{\rm train}= 2.95\cdot{10}^{-2}$ and $\mathcal L_{\rm test}=3.58\cdot 10^{-2}$, and a coefficient of determination of $R^2_{\rm train}=0.971$ and $R^2_{\rm test}=0.965$ \ADDED{for $n_s=19$, and losses of $\mathcal L_{\rm train} = 2.64\cdot 10^{-2}$, $\mathcal L_{\rm test} = 3.21\cdot 10^{-2}$ and coefficients of determination $R^2_{\rm train}=0.973$ and $R^2_{\rm test} = 0.968$ for $n_s=49$. Since both models perform similarly, we conclude that using $n_s=19$ represents sufficiently the geometric information in this DNN model without incurring in a large set of features. The maximum and minimum errors as defined in \eqref{eq:residual} are $(r_{\min}^2,r_{\max}^2)_{\rm train} = (1.05\cdot 10^{-10}, 7.84\cdot 10^{-1})$ and $(r_{\min}^2,r_{\max}^2)_{\rm test} = (2.96\cdot 10^{-10}, 7.30\cdot 10^{-1})$ for $n_s=19$.}The predictions are depicted in a scatter plot in \figurename~\ref{fig:results_dnn_example}.

Although the solution is (obviously) less precise than the toy example, a certain level of accuracy is reached. Thus, a surrogate model via a non FE-based approach is developed. By just scanning in three directions a given lattice to compute its effective areas, and weighing the specimen to obtain its volume, the equivalent Young's modulus can be effectively predicted.

\section{Graph Neural Networks model}
\label{sec:gnn}

Now we interpret the lattice as a graph $\mathcal G (\mathcal V,\mathcal E)$, where the joints are the nodes $\mathcal V$ (also referred to as vertices), and the struts are the edges $\mathcal E$. Features can be defined for both the nodes and the edges, and so the GNN model is built. This architecture is agnostic to the input data size, only requiring a graph and its associated features. As an illustration, a properly trained graph-level task GNN is able to make predictions for two lattices with different topologies i.e., regardless of the differences in the number of joints and/or struts.

\subsection{Features and architecture}

In contrast to the DNN model, we will only\ADDED{feed}the GNN\ADDED{with}geometrical features without physics or mechanics context, to assess whether the GNN is able to mimic the underlying equations, that is, equilibrium laws and computation of equivalent property. We define as node features the (three) spatial coordinates $(x_j,y_j,z_j)$ of the joint $j$, whereas the edge features are the coordinates of both endpoints i.e., $(x_{e0}, y_{e0}, z_{e0}), (x_{e1}, y_{e1}, z_{e1})$ -- where 0,1 represent both endpoints -- and the associated length $\ell_e$ of the bar/beam $e$ represented by such edge. This leads to a total of seven edge features. 

This way, we are providing structured, geometrical graph information to feed the GNN model. Since the base material used is the same (316 Stainless Steel), no material information is provided. Otherwise, it could be defined as a feature at the edge level.\ADDED{Note that heterogeneous lattices i.e., in which each member has its (potentially different) material defined, are admitted in these architectures, however, they are out of the scope of this paper.}Additionally, the variable to be predicted is the equivalent Young's modulus $E_{z}^{\rm eq}$ along $z-$axis. The dataset can be compactly represented as follows:

\[
\left\lbrace
\lbrace(x_j^{(i)}, y_j^{(i)}, z_j^{(i)})\rbrace_{j=1}^{n_j},
\lbrace(x_{e0}^{(i)}, y_{e0}^{(i)}, z_{e0}^{(i)}), (x_{e1}^{(i)}, y_{e1}^{(i)}, z_{e1}^{(i)}), \ell_e^{(i)}\rbrace_{e=1}^{n_e},
{E_{z}^{\rm eq}}^{(i)}
\right\rbrace_{i=1}^{N}.
\]

The architecture of the model comprises the following steps: (1) initial dense or fully connected (FC) layers applied separately for each node and edge, leading to subsequent node and edge embeddings, (2) two graph updates (or hops), updating the node embeddings, and (3) output layer based on pooling and FC layer. The details of these modules are displayed in \tablename~\ref{tab:gnn_architecture}.\ADDED{The hyper-parameters have been computed after running a shallow, brute-force hypertuning, varying the output shape of the different layers in a range within the order of magnitude of the nodes and edges features. The results obtained for similar configurations are rather close.}

\begin{table}[htb]
    \centering
    \begin{tabular}{c|cccc}
        \hline
         Module & Layer & Input shape & Output shape & Activation  \\ \hline
         Node embedding & Dense & $(3,)$ & $(5,)$ & PReLU \\
         \rowcolor{lightgray}\cellcolor{white}Edge embedding & Dense & $(7,)$ & $(5,)$ & PReLU \\
         & Message passing & $(15,)$ & $(10,)$ & selu \\
         Graph Update \#1 & Convolution (mean) & $(10\times n_m,)$ & $(10,)$ & $-$ \\
         & Dense & $(15,)$ & $(10,)$ & selu \\
         \rowcolor{lightgray}\cellcolor{white} & Message passing & $(25,)$ & $(10,)$ & selu \\
        \rowcolor{lightgray}\cellcolor{white}Graph Update \#2 & Convolution (mean) & $(10\times n_m,)$ & $(10,)$ & $-$ \\
        \rowcolor{lightgray}\cellcolor{white} & Dense & $(20,)$ & $(10,)$ & selu \\
        & Pooling (mean) & $(10\times n_j,)$ & $(10,)$ & $-$\\
         \multirow{-2}{*}{Output} & Dense & $(10,)$ & $(1,)$ & linear\\ \hline
    \end{tabular}
    \caption{Pin-jointed truss lattice GNN model. Architecture of the graph-level task regression, where the second component of the shapes is the batch size which the model is fed with. Note that $n_m$ in the convolution layer stands for the number of messages (adjacencies) arriving to one node in the graph -- different for each node --  and $n_j$ in the pooling layer is the number of joints (nodes) per graph---different for every graph.}
    \label{tab:gnn_architecture}
\end{table}

\ADDED{The input and output shapes displayed in \tablename~\ref{tab:gnn_architecture} operate as follows. There are 3 node features i.e., the coordinates, which are expanded into 5 through the node embedding (dense) model. Similarly, 7 edges features are compressed into 5 through the (dense) edge embedding. This way, both dense models present a balanced output shape---recall that the compression performed in the edge embedding does not lose critical information, since 6 out of 7 edge features are the coordinates, which have been already embedded in the nodes.}

\ADDED{In the first graph update, a message is passed through the edges using a dense layer which inputs the concatenation of the node embeddings of both endpoints (shapes: $5+5$) and the proper edge embedding (shape: 5), leading to a total of 15 features. The message passing (dense) layer has an output size of 10. Afterwards, a (mean) convolution is applied averaging the $n_m$ messages coinciding in all arrival nodes, thus performing $n_j$ convolutions. Then, an extra dense layer taking the old node state and the convoluted features (shape: $5+10$) as input, outputs the new node state, with each node presenting feature size of 10. The second graph update is performed analogously, arriving to new node states of size 10, likewise.}

\ADDED{Lastly, an element-wise pooling layer averages the features of all nodes in a graph, encapsulating the graph in a total of 10 features. Then, a final dense layer is applied to compress these 10 features into a unique output, which is the effective property predicted (per graph).}

\begin{figure}[H]
    \centering
    \includegraphics[width=\linewidth]{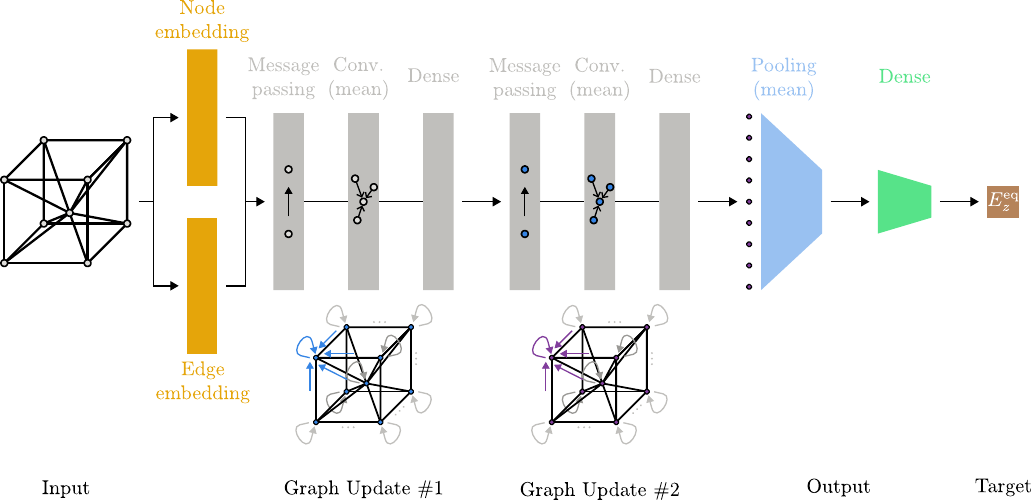}
    \caption{Pin-jointed truss lattice GNN model.\ADDED{Schematic}drawing of the graph model's architecture. The input is a lattice i.e., graph, and the goal is to predict an equivalent mechanical property of the lattice, namely the effective stiffness $E_z^{\rm eq}$ alongside $z-$axis. Initially, dense layers (depicted in yellow) are applied to every node and edge, leading to the subsequent node and edge embeddings. Then, the\ADDED{new}node features (or \ADDED{states)}are\ADDED{computed}with 2 graph updates (depicted in gray).\ADDED{Each of}these updates consists in a message passing where each edge computes a message by applying a dense layer to the concatenation of node states of both endpoints and the edge's own feature embedding; a convolution averaging at the common nodes of edges, and a dense layer to compute the new node state. This idea is illustrated in one arbitrary node of graph updates \#1 and \#2 using blue and purple arrows, respectively.\ADDED{Afterwards,}a pooling layer (blue) averaging the features of all nodes\ADDED{in a graph}is applied.\ADDED{In this example, the pooling layer takes the $n_j=9$ nodes of the structure, whose states are result of graph update \#2 (hence they are in purple). Lastly,} a dense layer (green)\ADDED{takes the result of this pooling and predicts}the output (brown),\ADDED{which is the effective property of the lattice i.e., graph.}}
    \label{fig:gnn_architecture}
\end{figure}

This architecture is sketched in \figurename~\ref{fig:gnn_architecture}. Dense layers (depicted in yellow) are initially applied to each node and edge in the graph, outputting node and edge embeddings of size 5. 

There are 2 graph updates (depicted in gray) which work as follows. (1) message passing: each edge computes a message by applying a dense layer to the concatenation of node states of both endpoints and the edge's own feature embedding, leading to an input shape of 15, and selecting 10 units as the output of the message passing. (2) convolution: messages are averaged (permutation-invariant pooling) at the common nodes of edges. (3) dense: at each node, a dense layer is applied to the concatenation of the old node state with the averaged edge inputs (15 in total) to compute the new node state, which is defined by 10 features. This idea is illustrated in the figure in one arbitrary node of graph updates \#1 and \#2 using blue and purple arrows, respectively. Then, the pooling layer (blue) averages the node features -- 10, given by the previous layer -- across all the nodes in a graph i.e., the whole graph is now encoded into 10 features. Lastly, a dense layer (green) is applied to get the unique output (brown) i.e., the equivalent mechanical property\ADDED{predicted.}

\subsection{Results for pin-jointed truss\ADDED{lattices}}

\begin{figure}[htb]
    \centering
    \includegraphics[width=.7\linewidth]{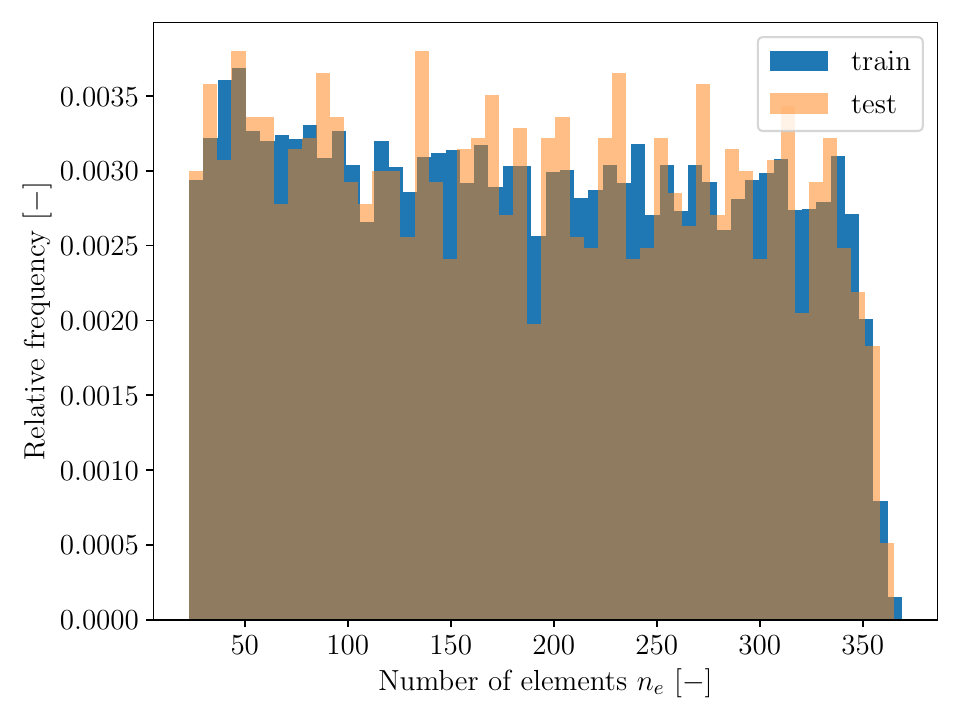}
    \caption{\ADDED{Pin-jointed truss lattices. Training and test set histograms of number of elements $n_e$ distribution.}}
    \label{fig:histograms}
\end{figure}

We first assess the performance of the model in pin-jointed truss lattices. The design of experiments (DoE) is performed similarly by generating lattices using the same material (316L Stainless Steel) and strut geometry (circular strut with radius $r$).\ADDED{The training set is generated by varying the number of non-fixed nodes from 1 to 50, sampling 400 lattices per each node, which leads to a total of $N=400\times 50 = 20\,000$ samples,}and 15\% of them are taken as validation set i.e., $3\,000$ samples.\ADDED{Separately, the test set is generated similarly, sampling 40 lattices per node, leading to 2\,000 samples. In this regard, the training and test sets have the same (uniform) distribution of number of nodes, whereas the number of elements $n_e$ varies according to \figurename~\ref{fig:histograms}.}

\ADDED{It can be observed in \figurename~\ref{fig:histograms} that the number of elements from the lattices in training and test sets are similarly distributed, indeed with slight differences due to the independent sampling processes.}

All the variables are scaled with standard normalization using the mean and standard deviation values of the training set. The optimization parameters are displayed in \tablename~\ref{tab:gnn_optimization_parameters}.

\begin{table}[htb]
    \centering
    \begin{tabular}{c|c}
         Parameter & Choice \\ \hline
         Optimizer & nadam \\
         \rowcolor{lightgray}Learning rate & $0.01$\\
         Epochs & 10\,000 \\
         \rowcolor{lightgray}Loss, $\mathcal L$ & MSE
    \end{tabular}
    \caption{Pin-jointed truss lattice GNN model. Optimization parameters of the training.}
    \label{tab:gnn_optimization_parameters}
\end{table}

\begin{figure}
    \centering
    \subfloat[Pin-jointed truss lattice GNN model.]
    {\includegraphics[width=.4\linewidth]{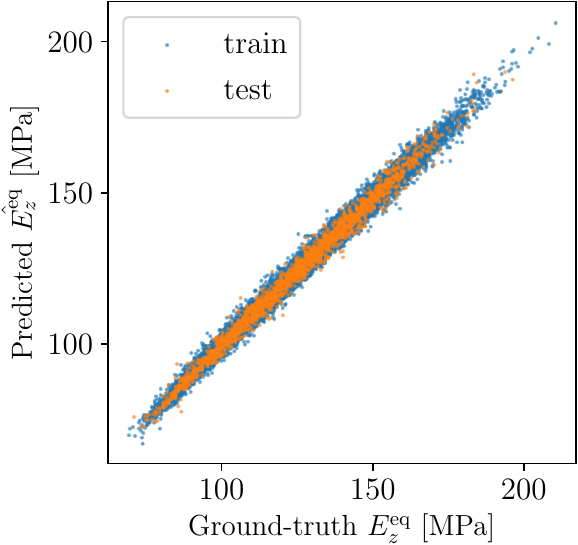}\label{fig:results_gnn_truss}}%
    \hspace{7.5ex}%
    \subfloat[Euler-Bernoulli beam lattice GNN model.]
    {\includegraphics[width=.4\linewidth]{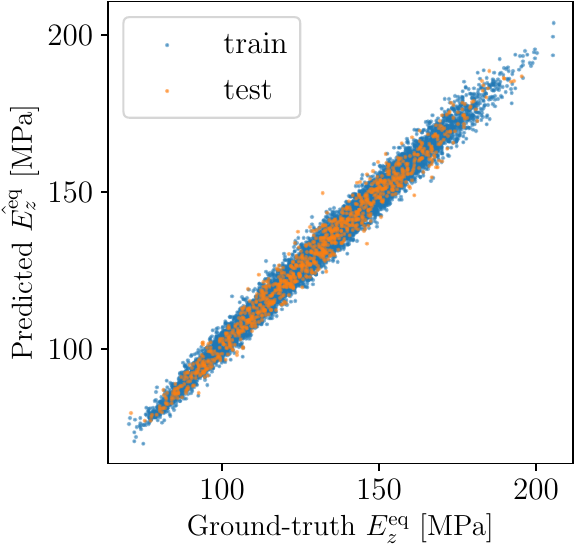}\label{fig:results_gnn_beam}}
    \caption{Pin-jointed truss (a) and Euler-Bernoulli beam (b) lattice GNN models. Equivalent stiffness $E_{z}^{\rm eq}$ prediction, where the train and test are represented by blue and orange points, respectively. Note that the diagonal line at $45^{\circ}$ is the perfect prediction.}
    \label{fig:results_gnn}
\end{figure}

The training yields $\mathcal L_{\rm train} = 1.02\cdot 10^{-2}$ and $\mathcal L_{\rm test} = 1.23\cdot 10^{-2}$, with coefficients of determination $R^2_{\rm train} = 0.9900$ and $R^2_{\rm test} = 0.9873$.\ADDED{The maximum and minimum errors as defined in \eqref{eq:residual} are $(r_{\min}^2,r_{\max}^2)_{\rm train} = (3.20\cdot 10^{-12}, 3.91\cdot 10^{-1})$ and $(r_{\min}^2,r_{\max}^2)_{\rm test} = (2.39\cdot 10^{-9}, 2.20\cdot 10^{-1})$.}The results of the prediction are depicted in \figurename~\ref{fig:results_gnn_truss}. Note the more accurate prediction in the test set than in the previous model, this being a sign that the GNN outperforms the DNN model.

\ADDED{To assess the learning capabilities of the model, we also run a KFold cross-validation consisting of $k=5$ folds performing a random split i.e., dropping 20\% of the information in five occasions. The results are displayed in \tablename~\ref{tab:kfold}. Note the similar performances of the 5 different training scenarios.}

\begin{table}[htb]
        \centering
        {\color{black}
        \begin{tabular}{c|ccccc|c}
                & $k=1$ & $k=2$ & $k=3$ & $k=4$ & $k=5$ & Whole dataset\\  \hline
             MSE$_{\rm test}$ & $1.47\cdot 10^{-2}$ & $1.48\cdot 10^{-2}$ & $9.66\cdot 10^{-3}$ & $1.91\cdot 10^{-2}$ & $1.51\cdot 10^{-2}$ & $1.23\cdot 10^{-2}$ \\
             $R^2_{\rm test}$ & $0.985$ & $0.985$ & $0.990$ & $0.981$ & $0.985$ & $0.987$
        \end{tabular}}
        \caption{\ADDED{Pin-jointed truss lattice GNN model. Results of the KFold cross-validation compared with the learning capabilities using the whole dataset.}}
        \label{tab:kfold}
\end{table}

Being all tests subjected to the same boundary conditions i.e., the uniaxial test to compute the equivalent property, the GNN is able to `solve' the equilibrium equations at the nodes by means of the initial node and edge embeddings and message-passing layers. Then, the output layer (pooling+dense) replicates the resultant of reaction forces -- reflected as $\bm f_r \cdot \bm u_r$ in Equation~\eqref{eq:equivalent_stiffness} -- which is a term needed to obtain the effective behavior of the lattice.

\subsection{Results for Euler-Bernoulli beam lattices}

Now we study a different mechanical model, namely a dataset of $N=10\,000$ Euler-Bernoulli beam lattices, with same material and geometrical properties than the pin-jointed truss lattice set counterpart. The same model as the previous section is applied, with same node and edge features, and architecture---see \tablename~\ref{tab:gnn_architecture} and \figurename~\ref{fig:gnn_architecture}. The validation set is 15\% of this dataset, whereas for the test set, $1\,000$ new samples (10\%) are generated separately. All features and targets are scaled with standard normalization, and the training is performed using the parameters displayed in \tablename~\ref{tab:gnn_optimization_parameters}. 

Once the GNN model is trained, the loss function evaluated at both datasets yields $\mathcal L_{\rm train} = 1.55\cdot 10^{-2}$ and $\mathcal L_{\rm test} = 2.20\cdot 10^{-2}$. Additionally, the coefficients of determination are $R^2_{\rm train} = 0.9845$ and $R^2_{\rm test} = 0.9778$.\ADDED{The maximum and minimum errors as defined in Equation \eqref{eq:residual} are $(r_{\min}^2,r_{\max}^2)_{\rm train} = (5.68\cdot 10^{-10}, 3.36\cdot 10^{-1})$ and $(r_{\min}^2,r_{\max}^2)_{\rm test} = (2.11\cdot 10^{-9}, 5.14\cdot 10^{-1})$.}The predictions are depicted in \figurename~\ref{fig:results_gnn_beam}. Although the prediction of this mechanical model becomes more challenging due to the appearance of shear forces and bending/twisting moments in the (rigid) joints -- leading to more equilibrium equations -- the GNN model is able to perform accurately on this dataset likewise. 

\ADDED{The observed performances are similar to benchmarks in the architected materials field e.g., Zheng et al. proposal to address it through variational autoencoders \cite{Zheng2023} ($R^2_{\rm test} = 0.982$); or the CNN architecture proposed by Shen et al. \cite{Shen2022} ($R^2_{\rm test} = 0.971$). A much simpler ML architecture (Random Forest) yielded particularly accurate results when training $6\,598$ periodic, fixed-topology, parameterized architectures (Euler-Bernoulli beams) $R^2_{\rm test} = 0.9991$ \cite{ben2023gam}. These are examples of the small error incurred by surrogates compared to the speed-up in computation time they allow.}

\section{Inverse problem}

Making use of the GNN surrogate models derived in Section~\ref{sec:gnn}, we now present an approach to address the inverse problem i.e., obtaining the lattice structure given a Young's modulus $E_z^{\rm eq, *}$ prescribed by the user. The model is a non-bijective function, thus the inverse problem is ill-posed \cite{lavrentiev1986ill,kaipio2006statistical}. Additionally, although a (potentially successful) Newton-Raphson approach to obtain numerically the inputs given an output may be derived -- the derivatives of the output w.r.t. the inputs can be computed via automatic differentiation \cite{baydin2018automatic} -- this approach is still not able to yield the number of nodes, struts and their adjacencies i.e., the graph $\mathcal G$---recall that the model inputs are the node and edge features, plus the adjacencies.

Therefore, thanks to the fast evaluations that the surrogate model allows, a database with $100\,000$ samples is rapidly generated to cover a wide range of Young's moduli $E_{z}^{\rm eq}$. With a sufficient number of samples predicted, several structural choices can yield the same effective property. That is the reason why we introduce the weight (or volume $V$ as defined in Equation~\eqref{eq:lattice_volume}, since the material stays the same) as an additional variable in this loop, which poses\ADDED{a}particular interest in structural optimization problems.

\begin{figure}[htb]
    \centering
    \includegraphics[width=.85\linewidth]{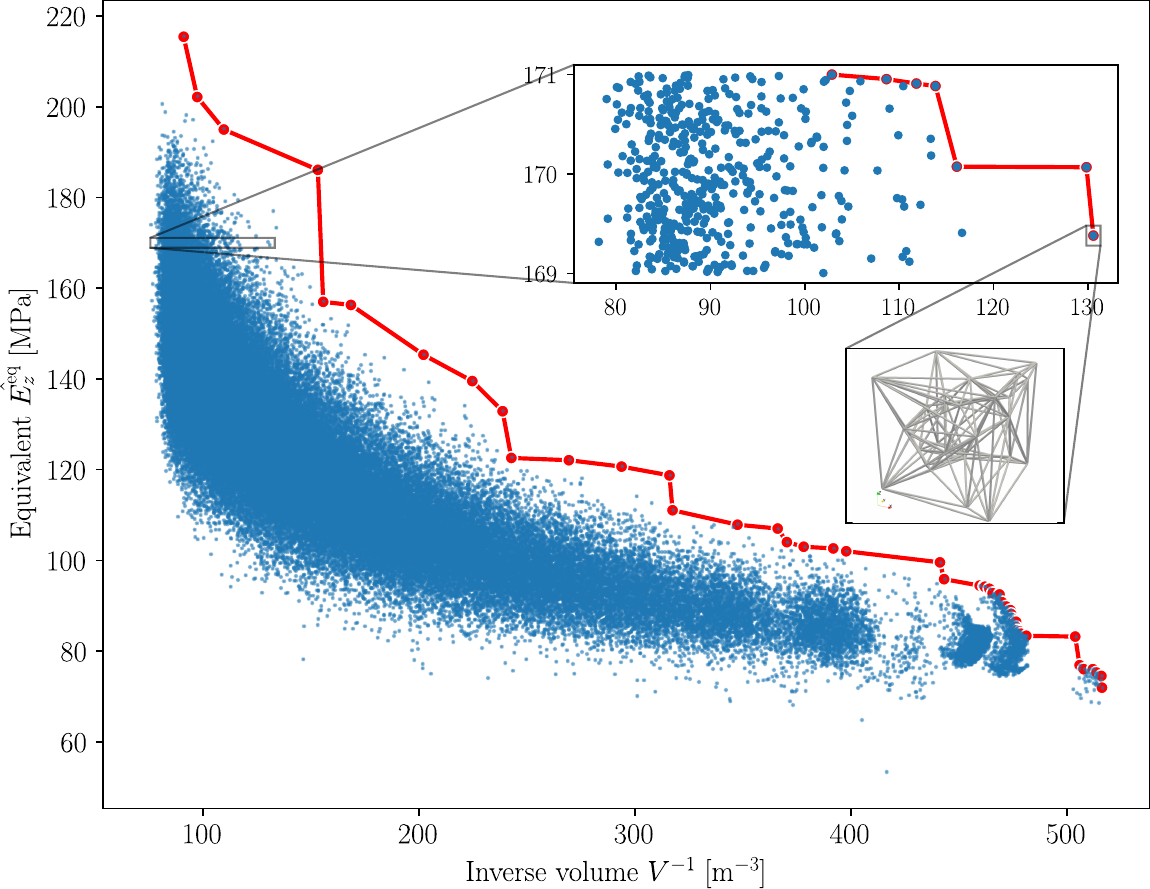}
    \caption{Distribution of inverse volume and equivalent predicted stiffness $(V^{-1},\hat{E_{z}^{\rm eq}})$ of $100\,000$ samples generated with the GNN surrogate model and its optimal Pareto front, in red. Local Pareto fronts can be defined likewise by introducing a range $[E_{\min}, E_{\max}]$, set to 169\,MPa and 171\,MPa in this case in the zoomed-in plot. Additionally, the structure with minimum weight of the local Pareto front is represented, with a predicted behavior of $\hat{E_z^{\rm eq}}=169.4\,$MPa, whereas the ground-truth property of such lattice yields $E_z^{\rm eq} = 171.0\,$MPa---implying an\ADDED{absolute value of the relative error $\epsilon$ with respect to the ground-truth property of $\epsilon=0.94\%$.}}
    \label{fig:pareto_front}
\end{figure}

Since lightweight structures are pursued, finding the minimum-weight lattice amongst all structures behaving with similar stiffness aims at this objective. To study the maximum-stiffness-minimum-weight configurations, every pair $\lbrace({{V}^{-1}}^{(i)}, {E_z^{\rm eq}}^{(i)})\rbrace_{i=1}^N$ is represented in a 2D scatter plot in \figurename~\ref{fig:pareto_front}. Then, a Pareto front \cite{boyd2004convex} is defined given the best stiffness-weight trade-offs, which is depicted in red in such figure.

When requesting a certain Young's modulus e.g., $E_{z}^{\rm eq,*} = 170\,$MPa, a local range $E_{z}^{\rm eq,*}\in\left[E_{\min},E_{\max}\right]$ is defined. Then, a local Pareto front is depicted to look for the lightest solution within that range. In the case of \figurename~\ref{fig:pareto_front}, $E_{\min}=169\,$MPa and $E_{\max}=171\,$MPa, thus the lightest solution corresponds to a lattice with $V=7.66\cdot 10^{-3}\,$m$^3$ and $\hat{E_z^{\rm eq}}=169.4\,$MPa. We note that the ground-truth effective behavior of such lattice given by the (direct) procedure described in Section~\ref{subsec:effective_behavior} yields $E_z^{\rm eq} = 171.0\,$MPa.\ADDED{By defining the absolute value of the relative error $\epsilon$ with respect to the ground-truth property as}
    \[
        \ADDED{\epsilon = \dfrac{|E_z^{\rm eq} - \hat{E_z^{\rm eq}}|}{E_z^{\rm eq}},}
    \]
\ADDED{the relative error for this specific case is $\epsilon = 0.94\%$.}

\ADDED{Similar inverse problem approaches via Variational Graph Auto-encoder yield $R^2>0.999$ \cite{Zheng2023} albeit with more geometric restrictions and a 10 times bigger dataset---$965\,736$ datapoints. An inverse design of aperiodic 2D structures by a Conditional Generative Adversarial Network (CGAN) \cite{Zheng2021} gives a MSE under $1\cdot 10^{-2}$ after 100 epochs with an equally big data set ($100\,000$ samples). An improved version of that CGAN generalized for 3D structures results in the same MSE after half the epochs (50), considering only $10\%$ of the dataset \cite{Zheng2023_3D}. Other GAN-powered lightweight optimization inverse processes show greater deviations i.e., $\mathrm{RMSE} > 0.2$ \cite{Challapalli2021}.}

\section{Conclusions}
We developed a surrogate model approach to fast compute effective mechanical properties in pin-jointed truss or Euler-Bernoulli beam (aperiodic) lattices with varied topologies. This work represents a framework to efficiently evaluate mechanical properties for many structures without incurring the prohibitive computational cost that classical FE approaches would take. Therefore, this framework is suitable for inverse analysis and, in turn, multi-scale structural optimization problems. The energetically equivalent effective behavior of Young's modulus is computed with DSM, yielding the ground-truth or high-fidelity solutions used to fit the surrogate model. Although the focus of this study is on predicting the equivalent Young's modulus, it may be straightforwardly generalized to other mechanical properties such as Poisson's ratio \cite{ben2023gam} or the homogenized components of the stiffness tensor \cite{kollmann2020deep}.

We first derive a DNN surrogate model, which slices the structure into $n_s$ portions in the three directions, considering the effective area $A_{i,s}^{\rm eff}$ per slice $s$ and direction $i$. We hypothesize it behaves similarly to a spring-in-series system, hence the slicing and computation of effective areas. Taking all the previous into account, with deep architectures presenting $\sim 10^5$ parameters, we can obtain sufficient accuracy for our surrogate model in terms of coefficient of determination $(R^{2}_{\rm test} = 0.965)$, successfully predicting the equivalent Young's modulus per unit volume $\overline{E_z^{\rm eq}}$. In other words, by just scanning in three directions a given lattice to compute its effective areas, and weighing the specimen to obtain its volume, the equivalent Young's modulus can be effectively predicted. 

Taking a step further, we build a GNN model to take advantage of the fact that lattices -- or any FE mesh -- are graphs and exploit the physical meaning of their adjacencies in the learning process. We embed a set of features in the nodes (coordinates) and edges (coordinates of endpoints and length) of the lattices built of the same material. The graph updates performed in the GNN consist of message-passing layers and posterior convolution, which replicate the equilibrium of forces at the nodes. These messages are sent through the elements (trusses, beams), accounting for geometrical (length, cross-sectional area) and material properties (Young's modulus). Lastly, the pooling layer accounts for the resultant reaction forces in the nodes required to compute the effective behavior $E_z^{\rm eq}$. These steps in the GNN are simple matrix operations -- all of them linear but the application of activation functions -- thus bypassing the classical FE analysis, which requires the assembly and inversion (or factorization) of the stiffness matrix to compute the displacements and reaction forces.

Regarding the results, a GNN with $\sim 10^3$ parameters can better predict the effective mechanical behavior $(R^2_{\rm test} = 0.987$ in case of pin-jointed truss lattices). Note that the GNN model requires two orders of magnitude less parameters than DNN model, achieving even better results in terms of coefficient of determination. We recall that the model is agnostic to the number of nodes and elements that a particular lattice structure presents, which represents a versatile feature \cite{sanchez2021gentle}.\ADDED{In the light of such remark and especially when comparing this approach to other similar studies \cite{Zheng2021,Challapalli2021,Zheng2023,Shen2022,Zheng2023_3D} commented throughout the result sections, performance is similar to the state-of-the-art, slightly better if we consider this study covers aperiodic, randomly generated 3D architectures---unlike the mentioned ones.}

The development of surrogate models allows for fast DoE generation aimed at covering a wide range of (equivalent) mechanical properties. Thus, if a homogeneous block of material with a specific stiffness is required -- as in topology optimization algorithms that output functionally graded metamaterials with variable stiffness, see \cite{ben2023topology, saucedo2023updated} -- the most compliant lattice might be found within this generated database. Therefore, an approach to solve inverse problems (macro-to-micro) is derived, making this framework suitable for inverse analysis in multi-scale structural optimization problems. Furthermore, it is possible to find several lattices fulfilling the exact stiffness requirement, leaving room to use another optimization criterion with more variables involved, e.g., the weight of the lattice. We develop global and local Pareto fronts to carry out this multi-objective optimization strategy. 

Although the whole study is restricted to linear elastic lattice structures of the same material, we note that this proof of concept represents a building block for considering non-linear, inelastic effects.\ADDED{For instance, path-independent non-linear settings e.g., hyperelasticity, would require to implement the strain energy curve of the member material as extra edge features of the graph. Likewise, an approach for path-dependent behavior e.g., plasticity, would also require quantifying the loading history, using the internal variables to this end. Another extension could be performed}for multi-material structures, even spanning different scales,\ADDED{since this architecture can address heterogeneous lattices i.e., with each member being made of (potentially different) materials.}This is of great practical relevance to address the data-driven description of deformable solids, with consequent valuable computational time savings.

\section*{Acknowledgements}
\noindent\begin{minipage}[c]{.15\linewidth}
\includegraphics[width=\linewidth]{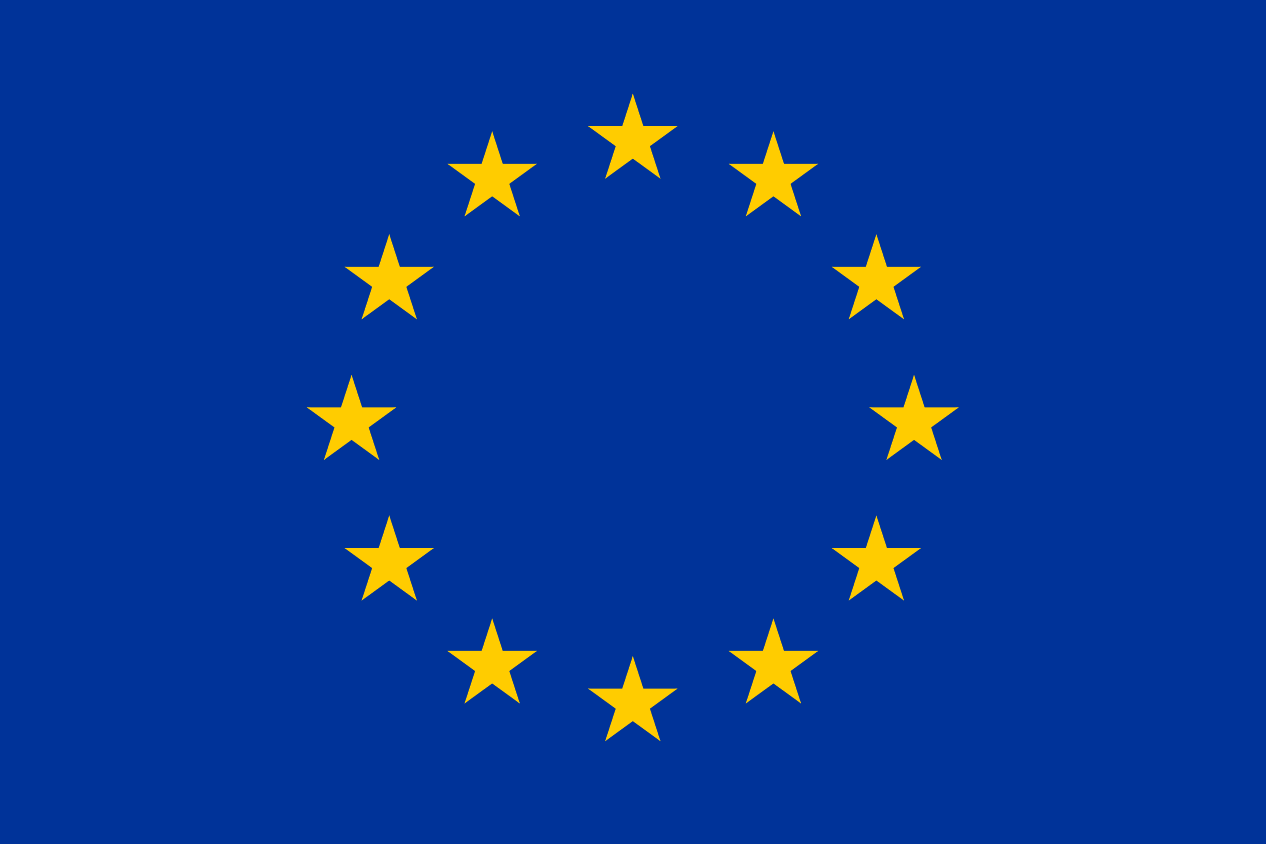}
\end{minipage}\hfill
\begin{minipage}[c]{.8\linewidth}
This project has received funding from the European Union's Horizon 2020 research and innovation programme under the Marie Skłodowska-Curie Grant Agreement No. 956401.
\end{minipage}

\section*{Declaration of interests}
The authors declare that they have no known competing financial interests or personal relationships that could have appeared to influence the work reported in this paper.

\bibliographystyle{unsrt}
\bibliography{gnn-mechanics}

\end{document}